March 22, 2021

# Convolutions for Stirling numbers, Lah numbers, and binomial coefficients

Khristo N. Boyadzhiev

*Department of Mathematics*
*Ohio Northern University*
*Ada, OH 45810, USA*

k-boyadzhiev@onu.edu

**Abstract**

A number of identities are proved mostly in the form of convolutions, involving Stirling numbers of both kinds, Lah numbers, and binomial coefficients. Certain "convolution polynomials" are discussed also. The proofs are based on series transformation formulas.

## 1. Introduction

Given two double sequences $a(n,k)$ and $b(n,k)$ where $n, k = 0, 1, 2, \ldots$, their convolution is defined by

$$c(n, p) = \sum_{k=0}^{n} a(n,k) b(k, p)$$

and represents a third double sequence. In general, this is not a commutative operation, we cannot interchange the two sequences.

In this article we carry out a systematic study of such convolutions for the double sequences of Stirling numbers of the first and second kinds, $s(n,k)$ and $S(n,k)$, the Lah numbers $L(n,k)$, and the binomial coefficients $\binom{n}{k}$. In this study we produce various identities relating the above numbers.

Excellent reference books for the Stirling and Lah numbers are the classical books of Comtet [15] and Graham et al. [17]. These numbers are very popular and important in mathematics. As additional references we will mention also [1, 3, 4, 5, 6, 8, 9, 10, 11, 12]. Certain convolution identities for the Stirling numbers of the first kind were considered by Agoh amd Dilcher in [2].

Our results are given in twelve propositions. Here is one sample, a convolution identity proved below in Proposition 10

$$\sum_{k=p}^{n} L(n,k)s(k,p) = (-1)^{n+p} s(n,p) \quad (n \geq p \geq 0).$$

We also study representations of some convolution polynomials

$$c(n,p)(\mu) = \sum_{k=0}^{n} a(n,k)b(k,p)\mu^k$$

and one of the results (Proposition 12) says that for any two integers $n \geq p \geq 0$ and any $0 < \mu < 1$

$$(-1)^{n-p} \sum_{k=p}^{n} s(n,k)S(k,p)\mu^k > 0.$$

The double sequences mentioned above have the following exponentials generating functions

(1) $$\frac{1}{p!}(e^x - 1)^p = \sum_{n=0}^{\infty} S(n,p) \frac{x^n}{n!}$$

(2) $$\frac{1}{p!}\ln^p(1+x) = \sum_{n=0}^{\infty} s(n,p) \frac{x^n}{n!}$$

(3) $$\frac{1}{p!}\left(\frac{x}{1-x}\right)^p = \sum_{n=0}^{\infty} L(n,p) \frac{x^n}{n!}$$

(4) $$\frac{x^p e^x}{p!} = \sum_{n=0}^{\infty} \binom{n}{p} \frac{x^n}{n!}.$$

Our main tools are the series transformation formulas listed below (see [3, 6, 15, 16]).

Given a function

(5) $$f(t) = \sum_{n=0}^{\infty} a_n \frac{t^n}{n!}$$

we have the following properties with appropriate parameters $\lambda, \mu$

(6) $$f\left(\frac{\mu}{\lambda}(e^{\lambda t} - 1)\right) = \sum_{n=0}^{\infty} \frac{t^n}{n!} \left\{ \sum_{k=0}^{n} S(n,k)\lambda^{n-k}\mu^k a_k \right\}$$

$$\text{(7)} \qquad f\left(\frac{\mu}{\lambda}\log(1+\lambda t)\right) = \sum_{n=0}^{\infty} \frac{t^n}{n!}\left\{\sum_{k=0}^{n} s(n,k)\lambda^{n-k}\mu^k a_k\right\}$$

$$\text{(8)} \qquad f\left(\frac{\mu t}{1-\lambda t}\right) = \sum_{n=0}^{\infty} \frac{t^n}{n!}\left\{\sum_{k=0}^{n} L(n,k)\lambda^{n-k}\mu^k a_k\right\}$$

$$\text{(9)} \qquad e^{\lambda t} f(t) = \sum_{n=0}^{\infty} \frac{t^n}{n!}\left\{\sum_{k=0}^{n} \binom{n}{k}\lambda^{n-k} a_k\right\}.$$

We will also use Euler's series transformation. Given an ordinary generating function

$$f(t) = \sum_{n=0}^{\infty} a_n t^n$$

we have

$$\text{(10)} \qquad \frac{1}{1-\lambda t} f\left(\frac{\mu t}{1-\lambda t}\right) = \sum_{n=0}^{\infty} t^n \left\{\sum_{k=0}^{n} \binom{n}{k}\lambda^{n-k}\mu^k a_k\right\}.$$

A simple rule for power series is also needed:

$$\text{(11)} \qquad \frac{1}{1-\lambda t}\sum_{n=0}^{\infty} a_n t^n = \sum_{n=0}^{\infty} \left\{\sum_{k=0}^{n} \lambda^{n-k} a_k\right\} t^n.$$

Our proofs are based on the formulas from (1) to (11). This is the analytical method.

Some convolutions for the sequences $s(n,k)$, $S(n,k)$, $L(n,k)$, and $\binom{n}{k}$ can be obtained (and have been obtained) by combinatorial methods. In this article we want to demonstrate that the analytical method is very efficient and universal.

We also state several open problems which hopefully will stimulate further research on such convolutions.

## 2. Convolutions for Stirling numbers and binomial coefficients

For worming up we start with two well-known identities.

**Proposition 1**. For any two nonnegative integers $n, p$ one has [17, p. 265]

$$\text{(12)} \qquad \sum_{k=0}^{n} \binom{n}{k} S(k,p) = S(n+1, p+1).$$

Proof. Applying property (9) with $\lambda = 1$ to the series (1) we have

$$\frac{1}{p!}e^x(e^x-1)^p = \frac{1}{p!}(e^x-1+1)(e^x-1)^p = \frac{p+1}{(p+1)!}(e^x-1)^{p+1} + \frac{1}{p!}(e^x-1)^p$$

$$= \sum_{n=0}^{\infty}(p+1)S(n,p+1)\frac{x^n}{n!} + \sum_{n=0}^{\infty}S(n,p)\frac{x^n}{n!} = \sum_{n=0}^{\infty}[(p+1)S(n,p+1)+S(n.p)]\frac{x^n}{n!}.$$

At the same time

$$\frac{1}{p!}e^x(e^x-1)^p = \sum_{n=0}^{\infty}\frac{x^n}{n!}\left\{\sum_{k=0}^{n}\binom{n}{k}S(k,p)\right\}$$

and comparing coefficients we get

$$\sum_{k=0}^{n}\binom{n}{k}S(k,p) = (p+1)S(n,p+1)+S(n.p) = S(n+1,p+1)$$

(the last equality is the familiar recurrence relation for the Stirling numbers of the second kind [17, p. 264, equation (6.15)]).

**Proposition 2**. For any two nonnegative integers $n, p$

(13) $\quad \sum_{k=0}^{n} s(n,k)(-1)^k \binom{k}{p} = (-1)^p s(n+1, p+1).$

This also can be found in [17, p. 265, equation (6.16)]. It is listed there in the form

$$\sum_{k=0}^{n}\begin{bmatrix}n\\k\end{bmatrix}\binom{k}{p} = \begin{bmatrix}n+1\\p+1\end{bmatrix}$$

where

$$\begin{bmatrix}n\\k\end{bmatrix} = (-1)^{n-k}s(n,k)$$

are the unsigned Stirling numbers of the first kind.

Proof. This time we use property (7) with $\mu=1, \lambda=-1$ for the generating function (4). This gives (using (11) in the last step)

$$\frac{(-1)^p}{p!}\ln^p(1-t)e^{-\ln(1-t)} = \frac{(-1)^p}{1-t}\frac{\ln^p(1-t)}{p!} = \frac{(-1)^p}{1-t}\sum_{n=0}^{\infty}(-1)^n s(n,p)\frac{t^n}{n!}$$

$$= \sum_{n=0}^{\infty}\frac{t^n}{n!}\left\{(-1)^p\sum_{k=0}^{n}(-1)^k\frac{s(k,p)n!}{k!}\right\}.$$

At the same time the left hand side equals

$$\sum_{n=0}^{\infty}\frac{t^n}{n!}\left\{\sum_{k=0}^{n}s(n,k)(-1)^{n-k}\binom{k}{p}\right\}.$$

Comparing coefficients gives

$$\sum_{k=0}^{n}s(n,k)(-1)^{n-k}\binom{k}{p}=n!\sum_{k=0}^{n}(-1)^{k-p}\frac{s(k,p)}{k!}=(-1)^{n-p}s(n+1,p+1)$$

(the last equality is another well-known property of the Stirling numbers of the first kind [17, p. 265, equation (6.21)]).

One good feature of this proof is that it also implies a similar identity which is not so well-known

**Proposition 3.** For any two integers $n, p \geq 1$

(14) $\quad \sum_{k=0}^{n}s(n,k)\binom{k}{p}=s(n,p)+ns(n-1,p).$

For the proof we use again (4) and (7), this time with $\lambda=\mu=1$

$$\frac{1}{p!}\ln^p(1+t)e^{\ln(1+t)}=(1+t)\frac{\ln^p(1+t)}{p!}=(1+t)\sum_{n=0}^{\infty}s(n,p)\frac{t^n}{n!}$$

$$=\sum_{n=0}^{\infty}s(n,p)\frac{t^n}{n!}+\sum_{n=0}^{\infty}(n+1)s(n,p)\frac{t^{n+1}}{(n+1)!}=\sum_{n=1}^{\infty}s(n,p)\frac{t^n}{n!}+\sum_{n=1}^{\infty}ns(n-1,p)\frac{t^n}{n!}$$

$$=\sum_{n=1}^{\infty}[s(n,p)+ns(n-1,p)]\frac{t^n}{n!}=\sum_{n=1}^{\infty}\frac{t^n}{n!}\left\{\sum_{k=0}^{n}s(n,k)\binom{k}{p}\right\}$$

(we can start the summation from $n=1$ because $s(0,p)=0$). All done!

Identity (14) was found by Alfred Schreiber [16] by a different method.

We can extend this identity by using (4) and (7) with $\lambda=1$ and arbitrary $\mu$. Thus

$$\frac{\mu^p}{p!}\log^p(1+t)e^{\mu\log(1+t)}=\frac{\mu^p}{p!}(1+t)^\mu\log^p(1+t)=\sum_{n=0}^{\infty}\frac{t^n}{n!}\left\{\sum_{k=0}^{n}s(n,k)\binom{k}{p}\mu^k\right\}$$

and expending the left hand side we have

$$\frac{\mu^p}{p!}(1+t)^\mu\log^p(1+t)=\mu^p\left\{\sum_{n=0}^{\infty}\binom{\mu}{n}t^n\right\}\left\{\sum_{n=0}^{\infty}s(n,p)\frac{t^n}{n!}\right\}$$

$$=\mu^p\sum_{n=0}^{\infty}t^n\left\{\sum_{k=0}^{\infty}\binom{\mu}{k}\frac{s(n-k,p)}{(n-k)!}\right\}=\mu^p\sum_{n=0}^{\infty}\frac{t^n}{n!}\left\{\sum_{k=0}^{\infty}\binom{n}{k}\binom{\mu}{k}s(n-k,p)k!\right\}.$$

Comparison of coefficients gives

(15) $\quad \sum_{k=0}^{n} s(n,k) \binom{k}{p} \mu^k = \mu^p \sum_{k=0}^{\infty} \binom{n}{k}\binom{\mu}{k} s(n-k, p) k!.$

Clearly (14) results from here immediately when $\mu = 1$ and (13) also follows after a simple computation. For $\mu = 2$ one gets

$$\sum_{k=0}^{n} s(n,k) \binom{k}{p} 2^k = 2^p [s(n, p) + 2n\, s(n-1, p) + n(n-1) s(n-2, p)].$$

Now we look for convolutions like (12) where in the second position we want to put the numbers $s(k, p)$. First we make one interesting observation.

**Observation**. Let for $|t| < 1$

$$f(t) = \frac{\ln^p(1+t)}{t} = \sum_{n=0}^{\infty} a_n t^n,\ f(-t) = \frac{\ln^p(1-t)}{-t} = \sum_{n=0}^{\infty} (-1)^n a_n t^n$$

Then

$$\frac{1}{1-t} f\left(\frac{t}{1-t}\right) = \frac{(-1)^{p-1} \ln^p(1-t)}{-t} = (-1)^{p-1} f(-t)$$

$$= (-1)^{p-1} \sum_{n=0}^{\infty} (-1)^n a_n t^n = \sum_{n=0}^{\infty} t^n \left[ \sum_{k=0}^{n} \binom{n}{k} a_k \right]$$

so that by comparing coefficients

(16) $\quad \sum_{k=0}^{n} \binom{n}{k} a_k = (-1)^{p-1}(-1)^n a_n$

Next we use the fact that

$$\frac{\ln(1+t)}{t} \sum_{n=0}^{\infty} c_n t^n = \sum_{n=0}^{\infty} t^n \left\{ \sum_{k=0}^{n} \frac{c_k (-1)^{n-k}}{n-k+1} \right\}$$

for any power series. This follows from (11) by integration with respect to $\lambda$ (see [6]).

From this

$$\ln(1+t) \sum_{n=0}^{\infty} a_n t^n = \sum_{n=0}^{\infty} t^{n+1} \left\{ \sum_{k=0}^{n} \frac{a_k (-1)^{n-k}}{n-k+1} \right\} = \sum_{m=1}^{\infty} t^m \left\{ (-1)^{m-1} \sum_{k=0}^{m-1} \frac{a_k (-1)^k}{m-k} \right\}.$$

Now we turn to the Stirling numbers of the first kind $s(n, p)$ with generating function (2) so that

$$\frac{\ln^p(1+t)}{t} = p!\sum_{n=p}^{\infty}\frac{s(n,p)}{n!}t^{n-1} = p!\sum_{n=0}^{\infty}\frac{s(n+1,p)}{(n+1)!}t^n.$$

**Proposition 4.** For any two nonnegative integers $n, p$ we have the convolution formula

(17) $$\sum_{k=0}^{n}\binom{n}{k}\frac{s(k+1,p)}{(k+1)!} = (-1)^{p-1}(-1)^n\frac{s(n+1,p)}{(n+1)!}.$$

For the proof we just use equation (16) with

$$a_k = \frac{s(k+1,p)}{(k+1)!}.$$

Next we write

$$\frac{\ln^{p+1}(1+t)}{t} = \ln(1+t)\frac{\ln^p(1+t)}{t} = \ln(1+t)p!\sum_{n=0}^{\infty}\frac{s(n+1,p)}{(n+1)!}t^n$$

$$= \sum_{n=1}^{\infty}t^n\left\{(-1)^{n-1}\sum_{k=0}^{n-1}\frac{(-1)^k}{n-k}p!\frac{s(k+1,p)}{(k+1)!}\right\}.$$

On the other hand

$$\frac{\ln^{p+1}(1+t)}{t} = (p+1)!\sum_{n=0}^{\infty}\frac{s(n+1,p+1)}{(n+1)!}t^n$$

so we get by comparing coefficients

$$(p+1)!\frac{s(n+1,p+1)}{(n+1)!} = (-1)^{n-1}\sum_{k=0}^{n-1}\frac{(-1)^k}{n-k}p!\frac{s(k+1,p)}{(k+1)!}$$

or

(18) $$(p+1)\frac{s(n+1,p+1)}{(n+1)!} = (-1)^{n-1}\sum_{k=0}^{n-1}\frac{(-1)^k}{n-k}\frac{s(k+1,p)}{(k+1)!}.$$

This can be written also in the form

(19) $$s(n+1,p+1) = \frac{(-1)^{n-1}(n+1)!}{p+1}\sum_{k=0}^{n-1}\frac{(-1)^k}{n-k}\frac{s(k+1,p)}{(k+1)!}$$

which compares nicely to the well-known property

$$s(n+1,p) = s(n,p-1) - n\,s(n,p).$$

Writing

$$\ln^{p+1}(1+t) = (p+1)!\sum_{n=0}^{\infty} \frac{s(n, p+1)}{n!} t^n$$

we find by differentiation

$$\frac{\ln^p(1+t)}{1+t} = p!\sum_{n=1}^{\infty} \frac{s(n, p+1)}{(n-1)!} t^{n-1} = p!\sum_{n=1}^{\infty} \frac{s(n+1, p+1)}{n!} t^n.$$

At the same time, from (11) we have

$$\frac{\ln^p(1+t)}{1+t} = p!\sum_{n=0}^{\infty} t^n \left\{\sum_{k=0}^{n}(-1)^{n-k} \frac{s(k, p)}{k!}\right\} = p!\sum_{n=0}^{\infty} \frac{t^n}{n!}\left\{n!\sum_{k=0}^{n}(-1)^{n-k} \frac{s(k, p)}{k!}\right\}$$

which gives by comparing coefficients

$$s(n+1, p+1) = n!\sum_{k=0}^{n}(-1)^{n-k} \frac{s(k, p)}{k!}.$$

This is the well-known property we used in the proof of Proposition 2.

We can write this property in the form of convolution

(20) $$s(n+1, p+1) = \sum_{k=0}^{n}\binom{n}{k}(-1)^{n-k} s(k, p)(n-k)! = \sum_{m=0}^{n}\binom{n}{m}(-1)^m s(n-m, p) m!.$$

Next we will show another convolution equation, an interesting counterpart to equation (13).

**Proposition 5**. For any two nonnegative integers $n, p$

(21) $$\sum_{k=0}^{n}\binom{n}{k}\frac{s(k, p)}{k!} = (-1)^p \sum_{k=0}^{n}\frac{(-1)^k s(k, p)}{k!} = \frac{(-1)^{n-p}}{n!} s(n+1, p+1).$$

Proof. Applying Euler's transformation to the generating function (2) we compute

$$\frac{1}{1-t}\ln^p\left(1+\frac{t}{1-t}\right) = \frac{(-1)^p \ln^p(1-t)}{1-t} = p!\sum_{n=0}^{\infty} t^n \left\{\sum_{k=0}^{n}\binom{n}{k}\frac{s(k, p)}{k!}\right\}.$$

At the same time

$$\ln^p(1-t) = p!\sum_{n=0}^{\infty} \frac{(-1)^n s(n, p)}{n!} t^n$$

and from (11)

$$\frac{(-1)^p \ln^p(1-t)}{1-t} = (-1)^p p!\sum_{n=0}^{\infty} t^n \left\{\sum_{k=0}^{n}\frac{(-1)^k s(k, p)}{k!}\right\}.$$

Comparing coefficients we come to equation (21).

Now we want to show why removing $k!$ from the left hand side in (21) is difficult.

Applying property (9), namely,

$$e^t \sum_{n=0}^{\infty} \frac{t^n}{n!} a_n = \sum_{n=0}^{\infty} \frac{t^n}{n!} \left\{ \sum_{k=0}^{n} \binom{n}{k} a_k \right\}.$$

to the generating function (2) we have

$$e^t \ln^p(1+t) = p! \sum_{n=0}^{\infty} \frac{t^n}{n!} \left\{ \sum_{k=0}^{n} \binom{n}{k} s(k,p) \right\}.$$

We have here a convolution without the factor $1/k!$, but the difficulty is to find the Taylor coefficients of the function $e^t \ln^p(1+t)$ in a different form. This represents an open problem – finding a concise form of the convolution

$$\sum_{k=0}^{n} \binom{n}{k} s(k,p).$$

Next we look for a convolution like (12) where $S(n,k)$ and the binomial change places. First we need to involve the exponential polynomials [9]

$$\varphi_n(x) = \sum_{k=0}^{n} S(n,k) x^k$$

with exponential generating function

$$e^{x(e^t - 1)} = \sum_{n=0}^{\infty} \varphi_n(x) \frac{t^n}{n!}.$$

**Proposition 6**. For any two non-negative integers $n \geq p$, and any $\mu$ we have

(22) $$\sum_{k=0}^{n} S(n,k) \binom{k}{p} \mu^k = \mu^p \sum_{k=0}^{n} \binom{n}{k} S(k,p) \varphi_{n-k}(\mu).$$

When $\mu = 1$ this becomes

(23) $$\sum_{k=0}^{n} S(n,k) \binom{k}{p} = \sum_{k=0}^{n} \binom{n}{k} S(k,p) b_{n-k}$$

where $b_n = \varphi_n(1)$ are the Bell numbers.

With $p = 0$ in (23) we have as expected

$$\sum_{k=0}^{n} S(n,k) = b_n.$$

When $p = 1$ in (23) we find

$$\sum_{k=0}^{n} S(n,k)k = \sum_{k=1}^{n} \binom{n}{k} b_{n-k} = \sum_{j=1}^{n-1} \binom{n}{j} b_j = b_{n+1} - b_n$$

which was obtained in [6]. The last equality follows from the property of the exponential polynomials

$$\varphi_{n+1}(x) = x \sum_{k=0}^{n} \binom{n}{k} \varphi_k(x)$$

(see [9]).

Proof of the proposition. We apply property (6) with $\lambda = 1$ to the generating function (4). Thus

$$\mu^p \frac{(e^t - 1)^p}{p!} e^{\mu(e^t - 1)} = \sum_{n=0}^{\infty} \frac{t^n}{n!} \left\{ \sum_{k=0}^{n} S(n,k) \binom{k}{p} \mu^k \right\}.$$

Expanding the left hand side we also have

$$\mu^p \frac{(e^t - 1)^p}{p!} e^{\mu(e^t - 1)} = \mu^p \left\{ \sum_{n=0}^{\infty} S(n,p) \frac{t^n}{n!} \right\} \left\{ \sum_{n=0}^{\infty} \varphi_n(\mu) \frac{t^n}{n!} \right\}$$

$$= \mu^p \sum_{n=0}^{\infty} \frac{t^n}{n!} \left\{ n! \sum_{k=0}^{n} \frac{S(n,k)}{k!} \frac{\varphi_{n-k}(\mu)}{(n-k)!} \right\} = \mu^p \sum_{n=0}^{\infty} \frac{t^n}{n!} \left\{ \sum_{k=0}^{n} \binom{n}{k} S(n,k) \varphi_{n-k}(\mu) \right\}$$

and the result follows by comparing coefficients.

Here is another open problem. After Proposition 6 a reasonable question appears: can we "put" $\mu^k$ in equation (12)?

We can write equation (1) with $\mu x$ in the place of $x$ and then apply (9) with $\lambda = 1$ to get

$$\frac{1}{p!} e^x (e^{\mu x} - 1)^p = \sum_{n=0}^{\infty} \frac{x^n}{n!} \left\{ \sum_{k=0}^{n} \binom{n}{k} S(k,p) \mu^k \right\}$$

but now a similar argument (as in the proof of (12)) does not seem feasible.

## 3. Convolutions for Lah numbers and binomial coefficients

Recall that the Lah numbers have the form

$$L(n,k) = \frac{n!}{k!}\binom{n-1}{k-1} = \left(\frac{n!}{k!}\right)^2 \frac{k}{n(n-k)!} \quad (n \geq k \geq 0)$$

where $L(n,0) = 0 \ (n > 0)$, $L(n,k) = 0 \ (n < k)$, $L(n,1) = n!$, $L(0,0) = 1$.

We also need the Laguerre polynomials $L_n^{(q)}(x)$ defined by the generating function

(24) $$\frac{1}{(1-t)^{q+1}} e^{-\frac{xt}{1-t}} = \sum_{n=0}^{\infty} L_n^{(q)}(x) t^n .$$

**Proposition 7**. For any two non-negative integers $n \geq p$ and any $\mu$ we have

(25) $$\sum_{k=p}^{n} L(n,k) \binom{k}{p} \mu^k = \frac{\mu^p n!}{p!} L_{n-p}^{(p-1)}(-\mu) .$$

When $p = 0$ this becomes

(26) $$\sum_{k=0}^{n} L(n,k) \mu^k = n! L_n^{(-1)}(-\mu)$$

which was found in [8].

Proof. We use property (8) with $\lambda = 1$, that is,

$$f\left(\frac{\mu t}{1-t}\right) = \sum_{n=0}^{\infty} \frac{t^n}{n!} \left\{ \sum_{k=0}^{n} L(n,k) \mu^k a_k \right\}$$

for the generating function (4). This gives

$$\frac{1}{p!} \left(\frac{\mu t}{1-t}\right)^p e^{\frac{\mu t}{1-t}} = \sum_{n=p}^{\infty} \frac{t^n}{n!} \left\{ \sum_{k=p}^{n} L(n,k) \mu^k \binom{k}{p} \right\} .$$

Now we represent the left hand side in a different way, recognizing there the generating function for the Laguerre polynomials (24)

$$\frac{1}{p!} \left(\frac{\mu t}{1-t}\right)^p e^{\frac{\mu t}{1-t}} = \frac{\mu^p t^p}{p!} \frac{1}{(1-t)^p} e^{\frac{\mu t}{1-t}} = \frac{\mu^p t^p}{p!} \sum_{m=0}^{\infty} t^m L_m^{(p-1)}(-\mu)$$

$$= \frac{\mu^p}{p!} \sum_{m=0}^{\infty} t^{m+p} L_m^{(p-1)}(-\mu) = \frac{\mu^p}{p!} \sum_{n=p}^{\infty} \frac{t^n}{n!} \left\{ n! L_{n-p}^{(p-1)}(-\mu) \right\}$$

(setting $m + p = n$ in the last step). The assertion follows by comparing coefficients.

An open problem: If we want to swap the positions of $L(n,k)$ and $\binom{n}{k}$ in (25) we need to use property (8) for the generating function (3). Thus

$$\frac{e^{\lambda x}}{p!}\left(\frac{x}{1-x}\right)^p = \sum_{n=p}^{\infty}\frac{x^n}{n!}\left\{\sum_{k=0}^{n}\binom{n}{k}L(k,p)\lambda^{n-k}\right\}$$

but here even for $\lambda = 1$ it is not clear how to find an appropriate expansion of the left hand side.

We will use property (10) for (3) instead. This brings to a very interesting identity.

**Proposition 8.** For any two integers $n \geq p \geq 0$

(27) $$\sum_{k=p}^{n}\binom{n}{k}L(k,p)\frac{(-1)^k}{k!} = \frac{(-1)^p}{p!}.$$

Notice that the sum does not depend on $n$. According to the definition of Lah numbers this equation can be written in the form

$$\sum_{k=p}^{n}\binom{n}{k}\binom{k-1}{p-1}(-1)^k = (-1)^p$$

which is a known binomial identity allowing also a combinatorial proof.

Proof. First we write (3) in the form (replacing $x$ by $-x$)

$$\frac{(-1)^p}{p!}\left(\frac{x}{1+x}\right)^p = \sum_{n=0}^{\infty}x^n\left\{\frac{(-1)^n}{n!}L(n,k)\right\}$$

and then we apply (10) where $\lambda = \mu = 1$. This gives

$$\frac{(-1)^p}{p!}\frac{t^p}{1-t} = \sum_{n=p}^{\infty}t^n\left\{\sum_{k=p}^{n}\binom{n}{k}\frac{(-1)^k}{k!}L(k,p)\right\}.$$

The left hand side we expand this way

$$\frac{(-1)^p}{p!}\frac{t^p}{1-t} = \frac{(-1)^p}{p!}\sum_{m=0}^{\infty}t^{p+m} = \frac{(-1)^p}{p!}\sum_{n=p}^{\infty}t^n$$

and comparing coefficients we complete the proof.

## 4. Convolutions for Stirling numbers and Lah numbers

A nice representation of Lah numbers in terms of Stirling numbers has been known for some time. It looks very neat when written with the "vertical" form of the Stirling numbers

(28) $$L(n, p) = \sum_{k=p}^{n} \begin{bmatrix} n \\ k \end{bmatrix} \begin{Bmatrix} k \\ p \end{Bmatrix} = \sum_{k=0}^{n} (-1)^{n-k} s(n,k) S(k, p) = (-1)^n \sum_{k=0}^{n} s(n,k) S(k, p)(-1)^k.$$

Using the inversion of the Stirling transform [6,11]

$$b_n = \sum_{k=0}^{n} S(n,k) a_k \quad \leftrightarrow \quad a_n = \sum_{k=0}^{n} s(n,k) b_k$$

we can write also

(29) $$\sum_{k=p}^{n} S(n,k) L(k, p)(-1)^k = (-1)^n S(n, p).$$

Next we present a different form of this convolution, which brings to a special representation for the Stirling numbers of the second kind.

**Proposition 9**. For any two integers $n \geq p \geq 0$

(30) $$\sum_{k=p}^{n} S(n,k) L(k, p)(-1)^k = (-1)^p \sum_{k=p}^{n} S(n,k) p^{n-k} (-1)^{n-k}.$$

and also

(31) $$S(n, p) = (-1)^p \sum_{k=p}^{n} S(n,k) p^{n-k} (-1)^k.$$

Here (31) follows from (29) and (30). We will prove now (29) and (30) by using (3) and (6). First we write (3) in the form

$$\frac{1}{p!} \left( \frac{-x}{1+x} \right)^p = \sum_{n=0}^{\infty} \frac{x^n}{n!} \left\{ (-1)^n L(n,k) \right\}$$

and then we apply property (6) to the function on the left hand side

$$\frac{1}{p!} \left( \frac{1-e^t}{e^t} \right)^p = \frac{1}{p!} (e^{-t} - 1) = \sum_{n=0}^{\infty} \frac{t^n}{n!} \left\{ (-1)^n S(n,k) \right\}$$

$$= \sum_{n=0}^{\infty} \frac{t^n}{n!} \left\{ \sum_{k=p}^{n} S(n,k) L(k, p)(-1)^k \right\}.$$

This proves (29). To prove (30) we write the left hand side differently

$$\frac{1}{p!} \left( \frac{1-e^t}{e^t} \right)^p = (-1)^p e^{-pt} \frac{(e^t - 1)^p}{p!} = (-1)^p e^{-pt} \sum_{n=0}^{\infty} \frac{t^n}{n!} S(n, p)$$

$$= (-1)^p \sum_{n=0}^{\infty} \frac{t^n}{n!} \left\{ \sum_{k=p}^{n} S(n,p)(-p)^{n-k} \right\}$$

(for the last equality we use (9)). This proves (30) by comparing coefficients.

The factor $(-1)^k$ cannot be easily removed from (30). The same approach gives

$$\frac{1}{p!}\left(\frac{e^t-1}{2-e^t}\right)^p = \sum_{n=0}^{\infty} \frac{t^n}{n!} \left\{ \sum_{k=p}^{n} S(n,k)L(k,p) \right\}$$

and it is not clear how to write the Taylor coefficients of the left hand side in a different form.

**Proposition 10**. For any two integers $n \geq p \geq 0$

(32) $\quad \sum_{k=p}^{n} L(n,k)s(k,p) = (-1)^{n+p} s(n,p)$.

For the proof we apply property (8) with $\lambda = \mu = 1$ to the generating function (2)

$$\frac{1}{p!}\ln^p\left(1+\frac{t}{1-t}\right) = \frac{(-1)^p}{p!}\ln^p(1-t) = (-1)^p \sum_{n=p}^{\infty} \frac{t^n}{n!}\{(-1)^n s(n,p)\}.$$

According to (8) the left hand side also has the form

$$\sum_{n=0}^{\infty} \frac{t^n}{n!} \left\{ \sum_{k=0}^{n} L(n,k)s(k,p) \right\}$$

and this is all we need to finish the proof.

Two open problems: Can we replace $s(k,p)$ by $S(k,p)$ in (32)? Only if we know how to deal with

$$\frac{1}{p!}\left(\exp\left(\frac{t}{1-t}\right)-1\right)^p = \sum_{n=0}^{\infty} \frac{t^n}{n!} \left\{ \sum_{k=0}^{n} L(n,k) S(k,p) \right\}.$$

Computing the convolution $\sum_{k=p}^{n} s(n,k) L(k,p)$ is also difficult. Using (3) in (7) gives

$$\frac{1}{p!}\left(\frac{\ln(1+t)}{1-\ln(1+t)}\right)^p = \sum_{n=0}^{\infty} \frac{t^n}{n!} \left\{ \sum_{k=p}^{n} s(n,k)L(k,p) \right\}$$

and it is not clear how to find another expansion of the left hand side. It is possible that these convolutions could be computed in a convenient form by other means.

## 5. Convolutions between Stirling numbers of the first and second kind

It is good to mention first the two classical convolutions

$$(33) \quad \sum_{k=0}^{n} S(n,k) s(k,j) = \begin{cases} 1 & n = j \\ 0 & n \neq j \end{cases} \quad \text{and} \quad \sum_{k=0}^{n} s(n,k) S(k,j) = \begin{cases} 1 & n = j \\ 0 & n \neq j \end{cases}.$$

The second formula has a nice extension found by Charalambides and Todorov [12, 13, 19].

**Proposition 11**. For any two integers $n \geq p \geq 0$ and every complex number $z$ we have the representation

$$(34) \quad \sum_{k=0}^{n} s(n,k) \, S(k,p) \, z^k = (-1)^p \frac{n!}{p!} \sum_{j=0}^{p} \binom{p}{j} (-1)^j \binom{zj}{n}.$$

A more recent proof based on a different method can be found in [10].

When $z = 1$ we have from (33)

$$\sum_{k=p}^{n} s(n,k) \, S(k,p) = (-1)^p \frac{n!}{p!} \sum_{j=n}^{p} \binom{p}{j} (-1)^j \binom{j}{n} = \begin{cases} 0 & p < n \\ 1 & p = n \end{cases}.$$

When $z = -1$ we have the equation (28) (see [8])

$$\sum_{k=0}^{n} s(n,k) S(k,p) (-1)^k = (-1)^n L(n,p).$$

Todorov [19] showed that

$$(35) \quad \frac{1}{p!}((1+t)^\mu - 1)^p = \sum_{n=k}^{\infty} \frac{t^n}{n!} \left\{ \sum_{k=p}^{n} s(n,k) S(k,p) \mu^k \right\}$$

by using (1) and (7). Then he also expanded the left hand side in (35) as binomial series to prove (34).

An open problem: It is desirable to find a similar expression for the polynomial

$$f_{n,p}(z) = \sum_{k=0}^{n} S(n,k) \, s(k,p) \, z^k.$$

Applying property (6) with $\lambda = 1$ to the generating function (2) we can write

$$\frac{1}{p!} \ln^p(1 + \mu(e^t - 1)) = \frac{1}{p!} \ln^p(1 - \mu + \mu e^t) = \sum_{n=0}^{\infty} \frac{t^n}{n!} \left\{ \sum_{k=0}^{n} S(n,k) s(k,p) \mu^k \right\}$$

which shows that the problem reduces to finding the Taylor coefficients of the left hand side in a different form.

We can rewrite (35) as

$$(36) \quad \frac{1}{p!}(1-(1-t)^{\mu})^p = \sum_{n=k}^{\infty} \frac{t^n}{n!}\left\{(-1)^{n-p}\sum_{k=p}^{n} s(n,k)S(k,p)\mu^k\right\}$$

and we can show that the Maclaurin coefficients of the function $(1-(1-t)^{\mu})^p$ are positive when $0 < \mu < 1$. Namely, we have the following proposition.

**Proposition 12.** For every $0 < \mu < 1$ and every $n \geq p \geq 0$ the coefficients in (36) are positive. That is,

$$(37) \quad (-1)^{n-p}\sum_{k=p}^{n} s(n,k)S(k,p)\mu^k > 0.$$

The restriction $\mu < 1$ is essential. For $\mu = 1$ and $n > p$ these sums are zeros in view of (33). Also, for $n = 4, p = 3, \mu = 3$ we have the counterexample

$$(-1)^{4-3}\sum_{k=p}^{n} s(4,k)S(k,3)3^k = -324.$$

*Proof.* Let $0 < \mu < 1$. We first consider the case $p = 1$ (the case $p = 0$ is trivial). Since

$$(1-t)^{\mu} = \sum_{n=0}^{\infty} \binom{\mu}{n}(-1)^n t^n$$

we can write

$$(38) \quad 1-(1-t)^{\mu} = \sum_{n=1}^{\infty} \binom{\mu}{n}(-1)^{n-1}t^n = \sum_{n=1}^{\infty} \frac{t^n}{n!}(-1)^{n-1}\mu(\mu-1)\dots(\mu-n+1)$$

$$= \sum_{n=1}^{\infty} \frac{t^n}{n!}\{\mu(1-\mu)(2-\mu)\dots(n-1-\mu)\}.$$

Obviously, the coefficients are positive. When we raise both sides in (38) to power $p$, the Taylor coefficients on the right hand side will again be positive, as products of positive numbers.

**Remark**. For every two integers $n > p > 0$ and every $\mu > 0$ we have

$$\sum_{k=p}^{n}(-1)^{n-k} s(n,k)S(k,p)\mu^k = \sum_{k=p}^{n} \begin{bmatrix} n \\ k \end{bmatrix} S(k,p)\mu^k > 0$$

because the unsigned Stirling numbers of the first kind $\begin{bmatrix} n \\ k \end{bmatrix}$ and the Stirling numbers of the second kind $S(k,p)$ are positive. Proposition 12 says that for $0 < \mu < 1$ we also have

$$(38) \quad \sum_{k=p}^{n}(-1)^{k+p}\begin{bmatrix}n\\k\end{bmatrix}S(k,p)\mu^{k} > 0.$$